\documentclass[p5]{elsarticle}

\usepackage{lineno}
\modulolinenumbers[5]
\usepackage{bm} 
\usepackage{lmodern}
\usepackage{moreverb}
\usepackage{color}
\usepackage{textcomp}
\usepackage{amsfonts,amssymb,amsbsy}
\usepackage{graphicx,empheq}
\usepackage{scrextend}

\newcommand{\y}{{ \bf y}}

\newcommand{\cLtwo}{{\mathcal L}_2 (0,L)}

\usepackage{pifont}
\usepackage[tableposition=top]{caption}
\def\mX{{\mathbb X}}
\def\Ltwo{{\mathbb L}^2 }

\newtheorem{rmk}{Remark}[section]
\newtheorem{thm}{Theorem}[section]

\newtheorem{prop}{Proposition}[section]

\newtheorem{lem}{Lemma}[section]

\newcommand{\mc}{\mathcal}
\newcommand{\mb}{\mathbf}

\def\e1{{\varepsilon_{1}}}
\def\b1{{\beta_{11}}}
\def\bp3{{\beta_{33}}}
\def\ep3{{\varepsilon_{3}}}

\MHInternalSyntaxOn
\providecommand*\phantomword[3][c]{%
\mathchoice
{\MT_phantom_word:NNnn #1\displaystyle {#2}{#3}}%
{\MT_phantom_word:NNnn #1\textstyle {#2}{#3}}%
{\MT_phantom_word:NNnn #1\scriptstyle {#2}{#3}}%
{\MT_phantom_word:NNnn #1\scriptscriptstyle {#2}{#3}}%
}
\def\MT_phantom_word:NNnn #1#2#3#4{%
\@begin@tempboxa\hbox{$\m@th#2#4$}%
\setlength\@tempdima{\widthof{$\m@th#2#3$}}%
\hbox{\hb@xt@\@tempdima{\csname bm@#1\endcsname}}%
\@end@tempboxa}
\MHInternalSyntaxOff

\begin{document}

\begin{frontmatter}

\title{Modeling and semigroup formulation of charge or current-controlled active constrained layer (ACL) beams; electrostatic, quasi-static, and fully-dynamic assumptions}
\author[Paestum]{Ahmet \"Ozkan \"Ozer}\ead{aozer@unr.edu}    

\address[Paestum]{Department of Mathematics \& Statistics, University of Nevada, Reno, NV 89557, USA}  

\begin{abstract}

 A three-layer active constrained layer (ACL) beam model, consisting of a piezoelectric elastic layer, a stiff layer, and a constrained viscoelastic layer, is obtained  for cantilevered boundary conditions by using the reduced Rao-Nakra sandwich beam assumptions through a consistent variational approach. The Rao-Nakra sandwich beam assumptions keeps the longitudinal and rotational inertia terms. We consider electrostatic, quasi-static and fully dynamic assumptions due to Maxwell's equations. For that reason, we first include all magnetic effects for the piezoelectric layer.  Two PDE models are obtained; one for the charge-controlled case and one for the current-controlled case.  These two cases are considered separately since the underlying control operators are very different in nature. For both cases, the semigroup formulations are presented, and the corresponding Cauchy problems are shown to be well- posed in the natural energy space.
\end{abstract}


\end{frontmatter}


\section{Introduction}
\begin{figure}
\begin{center}
\includegraphics[height=4cm]{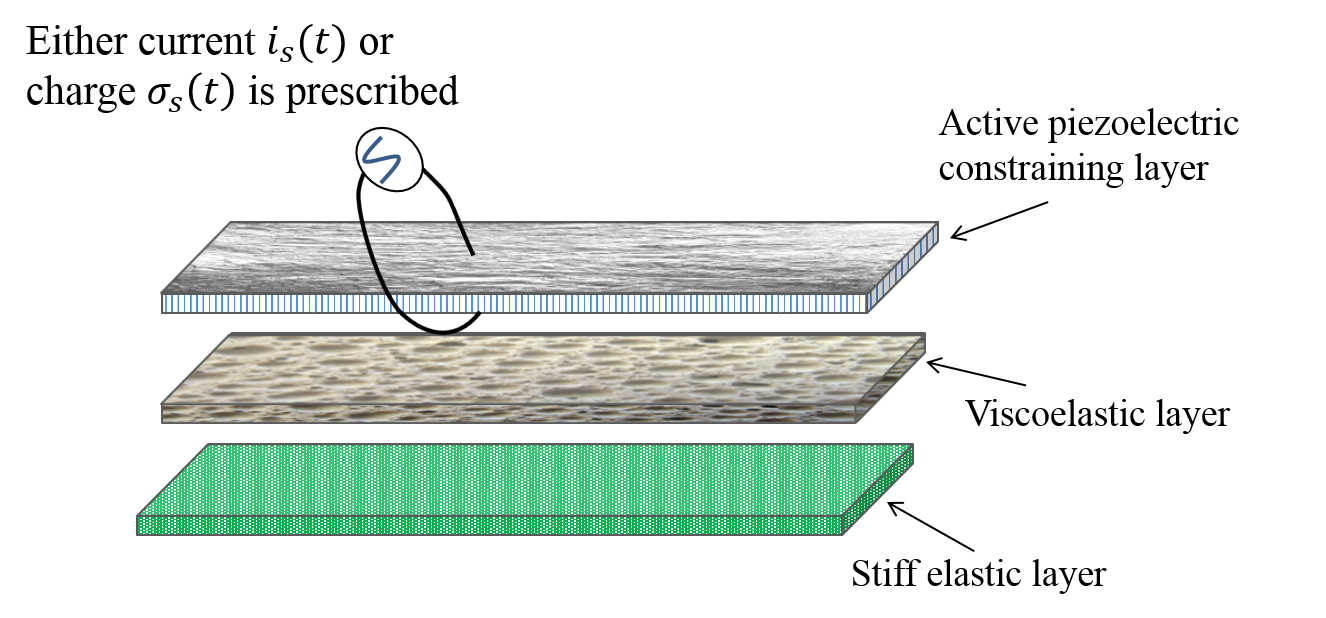}    
\caption{For a current or charge-actuated ACL beam, when   $i_s(t)$ or $\sigma_s(t)$ is supplied to the electrodes of the piezoelectric layer, an electric field is created between the electrodes, and therefore the piezoelectric beam either shrinks or extends, and this causes the whole composite stretch and bend.}  
\label{ACL}                                 
\end{center}                                 
\end{figure}

A three-layer actively constrained layer (ACl) beams is an elastic beam consisting of a stiff elastic beam and a piezoelectric beam constraining a viscoelastic beam which creates passive damping in the structure. The piezoelectric beam itself an elastic beam with electrodes at its top and bottom surfaces, insulated at the edges (to prevent fringing effects), and connected to an external
electric circuit (see Fig. \ref{ACL}).  These structures are widely used in  in civil, aeronautic and space space structures due to their small size
and high power density. They convert mechanical energy to electro-\emph{magnetic} energy, and vice versa. ACL composites involve piezoelectric layers and utilizes the benefits of them. Modeling these composites requires understanding the piezoelectric modeling since the ACL composites  are generally actuated through their piezoelectric layers. There are mainly  three ways to electrically actuate piezoelectric materials: voltage, current or charge. Piezoelectric materials have been traditionally activated by a voltage source \cite{Cao-Chen,Rogacheva,Smith}, and the references therein. However it has been observed that the type of actuation changes the controllability characteristics of the host structure. The main difference between three types of actuation  can be classified in three different groups:

\noindent \textbf{(i) Electrical nonlinearity (hysteresis):} Piezoelectric actuators actuated by a voltage source display nonlinear behavior at the higher drives. This is so-called hysteresis which results from the residual misalignment of crystal grains in the poled ceramic \cite{M-F}. For this reason, voltage-actuated piezoelectric beams are only actuated at low drives. This prevents them being used to their full potential. The hysteresis also causes a tracking error of 15\% of the displacement stroke, see \cite{Cao-Chen} the references therein.  It is observed that as the actuation is chosen charge or current instead of voltage, there is a substantial decrease in hysteresis, i.e. \cite{Cao-Chen,Review,Main1,Main2,M-F}, and the references therein.

\noindent \textbf{(ii) The bounded-ness of the control operator:}  Another difference is whether the underlying control operator in the state-space formulation of the control system is bounded  in the natural energy space. In the case of voltage and charge actuation, the underlying control operator is unbounded in the natural energy space, $\delta(\cdot)$ operators appear at the right side of the equations depending on the locations of sensors/actuators \cite{Hansen,O-M1}. However, the control operator is bounded in the case of current actuation, which is recently observed in \cite{O-M3} by using the potential formulation of the full set of Maxwell's equations in the modeling process.

\noindent  \textbf{(iii) Exact controllability \& Uniform stabilizability:}  In the case of a single piezoelectric beam, there is no way to exponentially stabilize the current controlled model with the magnetic effects \cite{O-M3}. Only strong stabilizability may be achieved for an infinitely many combinations of material parameters \cite{O-M1}. In the case of voltage actuation, the exponential \cite{O-M1} and polynomial \cite{Ozkan1}  stabilizability can be obtained for certain choices of parameters. This is the motivation behind \cite{O-M3} to derive charge or current-actuated distributed parameter models by using a thorough variational approach.

Accurately modeling an  ACL beam  also requires understanding of how  sandwich structures are modeled and how the interaction of the elastic layers are established. A three-layer sandwich beam consists of stiff outer layers and a  viscoelastic core layer. The core layer is supposed to deform only in transverse shear. The bending is uniform for the whole composite. Many sandwich beam models have been proposed in the literature, i.e., see \cite{Ditaranto,Mead,Rao,Sun,Trindade,Yan}, and the references therein.  These models mostly differ by the  assumptions on the viscoelastic layer. As well, depending on the inclusion of  the effects of longitudinal and rotational inertias, there are essentially two well-accepted models. The Mead-Marcus type models \cite{Mead} disregard these effects whereas the Rao-Nakra type models preserve them since it is  noted that these inertia effects are expected to have considerable importance especially at the high frequency modes for sandwich beams \cite{Rao}.

      The models for the ACL beams proposed in the literature mostly use the above-mentioned sandwich beam assumptions for the interactions of the layers, i.e. \cite{Stanway,Trindade}, and references therein.  The massive majority of these models are actuated by a voltage source, i.e. see \cite{Baz,F}, and the references therein. Moreover, the longitudinal vibrations were not taken into account in those papers; only the bending of the whole composite is studied.  In this paper, we use the reduced Rao-Nakra approach for which the longitudinal and rotational inertia terms are kept.  We even obtain a reduced Rao-Nakra model model by letting the weight and the stiffness of the middle layer go to zero \cite{Hansen3}.  Thus, a coupled but reduced PDE model is obtained for charge and current actuation. We also consider the assumptions of electrostatic, quasi-static and fully-dynamic in the Maxwell's equations. The biggest advantage of our models is being able to study controllability and stabilizability problems for ACL beams  in the infinite dimensional setting including all mechanical, electrical, and magnetic effects. We show that the proposed models with different types of actuation can be written in semigroup formulation, and they are shown to be well-posed in the natural energy space.


\section{Modeling Active-Constrained Layer
 (ACL) beams}

The ACL beam considered in this paper is a composite consisting of three layers that occupy the
region $\Omega=\Omega_{xy}\times (0, h)=[0,L]\times [-b,b] \times (0,h)$ at equilibrium where $\Omega_{xy}$ is a smooth bounded domain in the plane.
The total thickness $h$ is assumed to be small in comparison to the dimensions of $\Omega_{xy}$.
The beam  consists of a stiff layer, a compliant layer, and a piezoelectric layer, see Figure \ref{ACL}.

Let $0=z_0<z_1<z_2<z_3=h, $ with
 $$h_i=z_i-z_{i-1}, \quad i=1,2,3.$$
We use the rectangular coordinates $(x,y)$ to denote points in $\Omega_{xy},$ and  $(X, z)$ to denote points in $\Omega = \Omega^{\rm s} \cup \Omega^{\rm ve} \cup \Omega^{\rm p} $, where $\Omega^{\rm s}, \Omega^{\rm ve},$ and $\Omega^{\rm p}$ are the reference configurations of the stiff, viscoelastic, and piezoelectric layers, respectively, and they are defined by
\begin{eqnarray}
\nonumber &&\Omega^{\rm s}=\Omega_{xy}\times (z_0,z_1),\quad  \Omega^{\rm ve}=\Omega_{xy}\times (z_1,z_2), \quad  \Omega^{\rm p}=\Omega_{xy}\times (z_2,z_3).
\end{eqnarray}

 Define
$$\vec \psi=[\psi^1,\psi^2, \psi^3]^{\rm T}, \quad \vec \phi=[\phi^1,\phi^2,\phi^3]^{\rm T}, \quad \vec v=[v^1, v^2, v^3]^{\rm T}$$ where
\begin{eqnarray}
\label{defs1} &&  \psi^i=\frac{u^i-u^{i-1}}{h_i}, \quad \phi^i= \psi^i + w_x, \quad  v^i= \frac{u^{i-1}+u^i}{2}, \quad i = 1, 2, 3,\\
\label{defs3}&&\phi^1=\phi^3=0, \quad \quad \psi^1=\psi^3=-w_x, \quad \phi^2=\psi^2+ w_x
\end{eqnarray}

Let $T_{ij}$ and $ S_{ij}$   denote the stress and strain tensors for $i, j = 1, 2, 3$, respectively.
The constitutive equations for the piezoelectric layers are

\begin{eqnarray}
 \label{cons-eq50}
   T_{11}^{(3)}=\alpha^3 S^{(3)}_{11}-\gamma D_3, ~~   D_1=\varepsilon_{11} E_1, ~~  D_3=\gamma  S_{11}^{(3)}+\varepsilon_{33} E_3,
\end{eqnarray}
where $D, E, \alpha, \gamma,$  and $\varepsilon$ are electrical displacement vector, electric field intensity vector, elastic stiffness coefficient, piezoelectric coefficient,   permittivity coefficient, impermittivity coefficient, and and for the middle and the stiff layers are  given as the following:
\begin{eqnarray}
\label{cons-eq60}
 &  T_{11}^{(i)}=\alpha^i S_{11}^{(i)},\quad T_{13}^{(i)}= 2G_{2} S_{13}^{(i)}, \quad i=1,2&
\end{eqnarray}
where $G_2$ is the shear modulus of the second layer. Since we don't allow shear in the stiff layer we indeed have $T_{13}^{(i)}=0,$ $i=1,3.$ The strain components for the middle layer are
\begin{eqnarray}  \label{strains1} && S_{11}^{(2)}=\frac{\partial v^2}{\partial x}- (z-\hat z_i) \frac{\partial \psi^2}{\partial x}, \quad ~ S_{13}^{(2)}=\frac{1}{2}\left(\psi^2+ w_x\right)=\frac{1}{2}\phi^2,~
\end{eqnarray}
and for the piezoelectric and stiff layers are given by
\begin{eqnarray}
 \label{strains2} & S_{11}^{(i)}=\frac{\partial v^i}{\partial x}- (z-\hat z_i) \frac{\partial^2 w}{\partial x^2},\quad S^{(i)}_{13}=0,\quad i=1,3.&
\end{eqnarray}
\textbf{Magnetic effects and dynamic modeling:}      The full set of Maxwell's equations is (for instance see  \cite[Page 332]{Duvaut-L}):
\begin{eqnarray}
 \label{Maxwell} &&\left\{
  \begin{array}{ll}
 \nabla\cdot D =~\sigma_b  \quad{\rm{in}} \quad\Omega^{\rm p} \times \mathbb{R}^+~ & \text{(Electric Gauss's ~law)}\\
 \nabla\cdot B=~0 \quad{\rm{in}} \quad \Omega^{\rm p} \times \mathbb{R}^+~  &  \text{(Gauss's law of magnetism)}\\
 \nabla\times E=~-\dot B  \quad{\rm{in}} \quad \Omega^{\rm p} \times \mathbb{R}^+~ & \text{(Faraday's law)}\\
  \frac{1}{\mu}(\nabla\times B)= ~i_b +  \dot D  \quad{\rm{in}} \quad \Omega^{\rm p} \times \mathbb{R}^+~ & \text{(Amp\'{e}re-Maxwell law)}
  \end{array} \right.
\end{eqnarray}
with one of the essential electric boundary conditions prescribed on the electrodes of the piezoelectric layer
\begin{eqnarray}
 \label {vol} &&\left\{
  \begin{array}{ll}
 - D\cdot n =~\sigma_s   \quad{\rm{on}} \quad {\partial\Omega}^{\rm p} \times \mathbb{R}^+~\quad&\text{(Charge )}\\
 \frac{1}{\mu}(B \times n) =~i_s  \quad{\rm{on}} \quad {\partial\Omega}^{\rm p} \times \mathbb{R}^+~\quad & \text{(Current)}\\
 \phi=~V   \quad{\rm{on}} \quad {\partial\Omega}^{\rm p} \times \mathbb{R}^+~\quad& \text{(Voltage)}
  \end{array} \right.
\end{eqnarray}
and appropriate mechanical boundary conditions at the edges of the beam  (the beam is clamped, hinged, free, etc.).
Here $B$ denotes the magnetic field vector, and $\sigma_b, i_b, \sigma_s, i_s, V, \mu, n$ denote body charge density, body current density, surface charge density,
 surface current density, voltage, magnetic permeability, and unit normal vector to the surface $\partial\Omega^p,$ respectively. In this paper we consider  only current and charge-driven electrodes (i.e. we ignore the voltage boundary condition in (\ref{vol})). The voltage-driven electrode case is handled in details in \cite{Ozkan2}.

 By (\ref{Maxwell}), there exists a scalar electric potential $\varphi$ and a vector magnetic potential $A$ such that \begin{equation}
\label{imp1}B=\nabla\times A, \quad E=-\dot A-\nabla\varphi.
\end{equation}
 where $\dot A$ stands for the induced electric field due to the time-varying magnetic effects. In modeling piezoelectric beams, there are mainly three approaches including electric and magnetic effects \cite{Tiersten}:  Electrostatic, quasi-static, and fully dynamic electric field. The  electrostatic case assumes $B\equiv 0,$ and therefore $A\equiv 0.$
 The quasi-static  approach ignores $A$ and $\dot A$ since $A, \dot A\ll \varphi.$  With this assumption  $\dot D$ may be non-zero. However, in the fully dynamic case,   $A$ and $\dot A$ are left in the model.  Depending on the type of material, body charge density $\sigma_b$ and body current density $i_b$ can also be  non-zero. Note that even though the  displacement current $\dot D$ is assumed to be non-zero in both quasi-static and fully dynamic approaches, the term $\ddot D$  is  zero in quasi-static approach since $\dot A=0.$


 Henceforth, to simplify the notation, $x=x_1$ and $z=x_3.$

\noindent{\textbf{Electromagnetic assumptions.}} \label{Fully-dy}  The linear through-thickness assumption of the electric potential $\varphi(x,z)=\varphi^0(x)+z\varphi^1(x),$ which is a common assumption in many papers,  completely ignores the induced potential effect since $\varphi$ is completely known as a function of voltage. Therefore  we use a quadratic-through thickness potential distribution to account for this effect:
 \begin{eqnarray}
 \label{scalarpot} && \varphi(x,z)=\varphi^0(x)+z \varphi^1(x)+ \frac{z^2}{2} \varphi^2(x)
 \end{eqnarray}
Since we are in the beam theory, we assume that the magnetic vector potential $A$ has nonzero components only in $x$ and $z$ directions, i.e. $A_2\equiv 0.$ To keep the consistency with $\phi$, we assume that $A$ is quadratic through-thickness as well:
\begin{eqnarray}
  \label{vectorpot} && A_i(x,z)= A_i^0(x) + z A_i^1(x)  + \frac{z^2}{2} A_i^2(x), \quad i=1,3
\end{eqnarray}
By (\ref{imp1})
\begin{eqnarray}
\nonumber &&E_1= - \left( \dot A_1^0 +z \dot A_1^1 + \frac{z^2}{2} \dot A_1^2 \right)  -\left(({\varphi}^0)_x+z (\varphi^1)_x+ \frac{z^2}{2} (\varphi^2)_x\right),\\
\label{e3} &&E_3= -\left( \dot A_3^0 + z \dot A_3^1 +  \frac{z^2}{2} \dot A_3^2\right)-\left( \varphi^1 + z \varphi^2 \right).
\end{eqnarray}

\noindent {{\bf Work done by the external forces}} The work done by the electrical external forces (as in \cite{Lee,O-M3}) is
 {\footnotesize{
\begin{eqnarray}
\nonumber &&\mb{W}^{e}= \int_{\Omega^{\rm p}} \left(~-\tilde \sigma_b~ \phi+\tilde i_b \cdot A  \right) ~ dX + \int_{\partial \Omega^{\rm p}} \left(-\tilde\sigma_s~\phi + \tilde i_s\cdot  A\right)   ~ dX\\
\nonumber  &&= -\int_0^L\int_{-h_3/2}^{h_3/2} \tilde \sigma_b \left(\phi^0(x)+z \phi^1(x)+ \frac{z^2}{2} \phi^2(x)\right) ~dzdx\\
 \nonumber && + \int_0^L\int_{-h_3/2}^{h_3/2} \tilde i_b^1 \left(A_1^0(x) + z A_1^1(x)  + \frac{z^2}{2} A_1^2\right) ~ dz dx \\
 \nonumber && + \int_0^L  \left(-\tilde \sigma_s \left(\phi(h_3/2)-\phi(-h_3/2)\right) + \tilde i_s^1\left(A_1(h_3/2)-A_1(-h_3/2)\right)\right)  ~dx \\
\label{work}  &&= \int_0^L  \left(  -\sigma_b\left(\phi^0 + \frac{h_3^2}{24}\phi^2\right)-\sigma_s+i_b^1\left(A_1^0 + \frac{h_3^2}{24}A_1^2\right)- \sigma_s ~\phi^1 +i_s^1 ~A_1^1 \right)~ dx
\end{eqnarray}
}} where  $i_s(x,t)=(i_s^1(x,t),0,0),$ and $i_b(x,t)=(i_b^1(x,t), 0, i_b^3(x,t)). $  In the above $i_s$ has only one nonzero component since  $i_s \perp B $, and $i_s \perp n$ by (\ref{vol}). Moreover, $i_b$ has only one nonzero component since we assumed that there is no force acting in the $x_2$ and $x_3$ directions. Several remarks are in order:

 \begin{rmk} (i) We choose either surface charge $\sigma_s $ or $i_s$ to be non-zero, or $i_b$ or $\sigma_b$ to be nonzero, depending on the type of actuation. For consistency, we keep them both in deriving PDEs for the rest of the paper; in the case of only surface current actuation $i_s\ne 0, \sigma_s=i_b=\sigma\equiv 0,$ and in the case of surface charge actuation $\sigma_s\ne 0, i_s=i_b=\sigma_b\equiv 0.$\\

 (ii) Since the piezoelectric materials are not perfectly  insulated, the electric field $E$ causes currents to flow when conductivity occurs. Therefore the time-dependent equation of the continuity
of electric charge must be employed:
\begin{eqnarray}\nonumber &&\dot \sigma_b+ \nabla \cdot i_b =0 \quad {\text in} ~~\Omega^{\rm p}\\
\label{elec-cont2} &&\dot \sigma_s - i_b^3= 0, \quad {\rm or,}\quad\frac{ di_s}{dx}- i_b^3= 0 \quad \text{on}~ \partial\Omega^{\rm p}.
\end{eqnarray}
For more details, the reader can refer to \cite[Section 3.9]{Eringen}.\\

(iii) Note that if the magnetic effects are neglected, a variational approach cannot be used since $A\equiv 0$ and so the surface current may not be involved in $\mb W^e.$
This is very different from the charge and voltage actuation cases since the charge and voltage terms are not affected by $A$ in $\mb W^e.$
 \end{rmk}

Assume that the beam is subject to a  distribution of forces $(\tilde g^1, \tilde g^3, \tilde g)$ along its edge $x=L.$ In parallel to \cite{Hansen3}, then the total work done by the mechanical external forces is
\begin{eqnarray*}
\label{work-done}  \mb{W}^m&=&   g^1 v^1 (L)  + g^3 v^3 (L) +  g w(L)-M w_x(L).
\end{eqnarray*}

Now we use the constitutive equations (\ref{cons-eq50})- (\ref{strains2}), and (\ref{scalarpot})-(\ref{e3}) to find 
\begin{eqnarray}
\nonumber &&\mb P=\mb{P}^1 + \mb{P}^2=\frac{1}{2}\int_{\Omega^{\rm s} \cup \Omega^{\rm ve}} \sum_{i=1,2}\left(T_{11}^{(i)}S_{11}^{(i)}+ T_{13}^{(i)}S_{13}^{(i)}\right)~dX \\
\nonumber &&= \frac{1}{2} \int_0^L \left[{\alpha}^1 h_1\left((v^1_x)^2+ \frac{h_1^2}{12}w_{xx}^2\right)+{\alpha}^2 h_2\left((v^2_x)^2+ \frac{h_2^2}{12}(\psi_x)^2\right)+ G_2h_2 (\phi^2)^2\right]~dx,
\end{eqnarray}
\begin{eqnarray}
\nonumber &&\mb{E}^3-\mb{P}^3=\frac{1}{2}\int_{\Omega^{\rm p}} \left(D_1 E_1+D_3 E_3- T_{11}^{(3)}S_{11}^{(3)}-T_{13}^{(3)}S_{13}^{(3)}\right)~dX \\
\nonumber &=& \frac{1}{2} \int_0^L \left[-{\alpha}^3 h_3\left((v^3_x)^2+ \frac{h_3^2}{12}(w_{xx})^2\right)\right.\\
\nonumber && - 2{\gamma}h_3  \left(\left(\varphi^1+ \dot A_3^0 + \frac{h_3^2}{24} \dot A_3^2\right)v^3_x-\frac{h_3^2}{12}w_{xx}\left(\varphi^2+\dot A_3^1\right)\right) \\
\nonumber &&  + {\e1}h_3 \left( (\varphi^0_x)^2+\frac{h_3^2}{12} (\varphi^1_x)^2+ \frac{h_3^4}{320} (\varphi^2_x)^2 +(\dot A_1^0)^2+ \frac{h_3^2}{12} (\dot A_1^1)^2+ \frac{h_3^4}{320} (\dot A_1^3)^2  \right)\\
\nonumber &&   + {\ep3}h_3 \left( (\varphi^1)^2+\frac{h_3^2}{12} (\varphi^2)^2 +(\dot A_3^0)^2+ \frac{h_3^2}{12} (\dot A_3^1)^2+ \frac{h_3^4}{320} (\dot A_3^2)^2  \right)\\
\nonumber &&  + 2{\e1}h_3 \left( (\varphi^0)_x \dot A_1^0 + \frac{h_3^2}{24} (\varphi^0)_x( \varphi^2)_x + \frac{h_3^2}{24} \dot A_1^0 \dot A_1^2+ \frac{h_3^2}{24} (\varphi^0)_x \dot A_1^2 \right.\\
\nonumber &&\left. + \frac{h_3^2}{24} (\varphi^2)_x \dot A_1^0 +   \frac{h_3^2}{12} (\varphi^1)_x \dot A_1^1+  \frac{h_3^4}{320} (\varphi^2)_x \dot A_1^2  \right)\\
\label{E-P} &&\left.  + 2{\ep3}h_3 \left( \varphi^1 \dot A_3^0 + \frac{h_3^2}{24} \dot A_3^0 \dot A_3^2 + \frac{h_3^2}{24} \varphi^1  \dot A_3^2 + \frac{h_3^2}{12} \varphi^2 \dot A_3^0 \right)\right]~dx,
\end{eqnarray}
and $\mb K=\mb K^1+\mb K^2+\mb K^3 $ where
 \begin{eqnarray}
\label{KK1} \mb{K}^i&=&  \frac{\rho_i h_i}{2} \int_0^L \left((\dot v^i)^2 + \dot w^2 +\frac{h_i^2}{12}\dot w_x^2\right)~dx,\quad i=1,3,\\
\label{KK2} \mb{K}^2&=&  \frac{\rho_2 h_2}{2} \int_0^L \left((\dot v^2)^2 + (\dot \psi^2)^2+ \dot w^2 \right)~dx,
\end{eqnarray}
Finally, the magnetic energy of the piezoelectric beam is
\begin{eqnarray}
\nonumber &&\mb{B}^3=\frac{\mu}{2}\int_{\Omega^{\rm p}} (\nabla \times A) \cdot  (\nabla \times A) ~dX \\
\nonumber && =\frac{\mu h_3}{2} \int_0^L\left[ (A_1^1)^2+ \frac{(h_3)^2}{12} (A_1^2)^2 +((A_3^0)_x)^2 +\frac{(h_3)^2}{12}((A_3^1)_x)^2+ \frac{(h_3)^4}{320}((A_3^2)_x)^2 \right. \\
\label{BB} && \left.-2\left( A_1^1 ~(A_3^0)_x +\frac{(h_3)^2}{24} A_1^1(A_3^2)_x - \frac{(h_3)^2}{12} A_1^2(A_3^1)_x -\frac{(h_3)^2}{24}(A_3^0)_x(A_3^2)_x\right)\right]~dx.
\end{eqnarray}

\subsection{Variational Principle \& Equations of Motion}
By using (\ref{defs1})-(\ref{defs3}), the variables $\{v^2, \phi^2, \psi^2\}$ can be written as the functions of $\{v^1, v^3, p\}.$
As well, we know that the variables $\{\varphi^0, \varphi^2, A_3^1, A_1^0, A_1^2\}$ have nothing to do with the stretching equations as the single piezoelectric beam is actuated through the charge or the  current source \cite{O-M3}. For that reason, we choose only $\{v^1, v^3, w,   \phi^1, A_1^1,  A_3^0, A_3^2\}$ as the state variables.

To model charge or current-actuated ACL beams with magnetic effects, we use the following Lagrangian \cite{Lee,O-M3}
\begin{eqnarray}\label{Lag}  \mb{ L}= \int_0^T \left[\mb{K}-(\mb{P}-\mb{E}+\mb B)+ \mb{W}\right]dt\end{eqnarray}
where $\mb P-\mb E+\mb B=\mb P^1 + \mb P^2 +\mb P^3 -\mb E^3 + \mb B^3 $ is called electrical enthalpy where $\mb B^1=\mb B^2=\mb E^1=\mb E^2\equiv 0,$ and $\mb W=\mb W^m+\mb W^e$ is the total work done by the mechanical and electrical external forces, respectively. Note that in modeling piezoelectric beams by voltage-actuated electrodes we use a modified Lagrangian  \cite{O-M1}. Let $H=\frac{h_1 + 2h_2+h_3}{2}.$ We assume that the piezoelectric beam is clamped at $x=0$ and free at $x=L$.

The application of Hamilton's principle, setting the variation of Lagrangian $\mb L$ in  (\ref{Lag}) with respect to the all kinematically  admissible displacements of the chosen state variables to zero, yields a strongly coupled equations for the longitudinal and transverse dynamics together with the magnetic and electrical equations. It is not easy to study the controllability/stabilizability  properties. For this reason we study the following reduced model.

\subsection{Gauge condition and Rao-Nakra model assumptions}
 The magnetic potential vector $A$ and the electric potential $\phi$ are not uniquely defined in (\ref{imp1}).
  In fact, the Lagrangian $\mb L$  (\ref{Lag}) is invariant under a large class of transformations.
 \begin{thm}
 For any scalar $C^1$ function $\chi=\chi(x, z, t),$  the Lagrangian $\mb L$ is invariant invariant under the  transformation
\begin{eqnarray}
\label{Gaugemap}  & A \mapsto \tilde A:= A+\nabla \chi,\quad\quad \phi \mapsto \tilde \phi:= \phi-\dot\chi.&
\end{eqnarray}
\end{thm}

\textbf{Proof}: See \cite{O-M3}.

An additional condition can be added to remove the ambiguity as in \cite{O-M3}. The additional condition is generally known as  a {\em gauge} and it can be chosen to decouple the electrical potential equation  from the equations of the magnetic potential.
Define \begin{eqnarray}
\label{eta-theta}\eta:= A_3^0 + \frac{h_3^2}{24} A_3^2, \quad \theta:=A_1^1, \quad \xi:=\frac{\e1 h_3^2}{12\ep3}.
\end{eqnarray}
In conjunction to \cite{O-M3}, the appropriate gauge condition to decouple the magnetic and electric equations  is chosen
\begin{eqnarray}
\label{Colombo-trans}& -\xi \theta_x + \eta= 0, \quad \theta(0)=\theta (L)= 0 .
\end{eqnarray}

In this section we derive the  Rao-Nakra type ACL beam model using the  three-layer version of the Rao-Nakra sandwich beam model developed in \cite{Hansen3} by letting $\rho_2, \alpha^2\to 0.$ This
approximation retains the potential energy of shear and transverse kinetic energy. This together with the gauge condition (\ref{Colombo-trans}), we obtain the following model

\begin{eqnarray}
\label{reduced-nog} \left\{
  \begin{array}{ll}
   m \ddot w - K_1 \ddot{w}_{xx}+K_2w_{xxxx} -  H G_2   \phi^2_x=0&\\
  \rho_1h_1 \ddot v^1  -\alpha^1 h_1 v^1_{xx}  - G_2  \phi^2 = 0,   & \\
  \rho_3h_3 \ddot v^3       -\alpha^3 h_3 v^3_{xx}  + G_2 \phi^2 -{\gamma}h_3 \left((\varphi^1)_x +  \dot \eta_x \right) = 0,   & \\

    - \xi(\phi^1)_{xx} +\phi^1 -\frac{{\gamma}}{\ep3} v^3_x= \frac{\sigma_s (t) }{h_3\ep3} & \\

    \xi\ep3h_3\ddot \theta +\mu h_3 \left(\theta -\eta_x\right)+ \xi\ep3 h_3 (\dot \phi^1)_x=i_s^1(t) & \\

   \ep3 h_3 \ddot\eta + \mu h_3\left(\theta-\eta_x\right)_x + \ep3 h_3\dot \phi^1-\gamma h_3 \dot v^3_x=0&\\
      \phi^2 =\frac{1}{h_2}\left(-v^1+v^3\right) + \frac{H}{h_2}w_x.&
  \end{array} \right.
  \end{eqnarray}

with the natural boundary  conditions at $x=0,L$
\begin{eqnarray}
&&  \begin{array}{ll}
\nonumber   \left\{v^1, v^3, w, w_x \right\}_{x=0}=0, ~~   \left\{\alpha^1 h_1v^1_x\right\}_{x=L}=g^1,&\\
  \end{array} \\
  &&  \begin{array}{ll}
\nonumber  \left\{ \alpha^3 h_3 v^3_x  +{\gamma}h_3\left(\phi^1+\dot\eta\right)\right\}_{x=L}=g^3,~~ \left\{  \phi^1_x,~ \theta,~  \eta_x \right\}_{x=0,L}=0, &\\
  \end{array} \\
      &&  \begin{array}{ll}
\label{reduced-BC-nog}     \left\{K_2 w_{xx}  \right\}_{x=L} =-M, ~ \left\{K_1 \ddot w -K_2 w_{xxx} + G_2 H \phi^2 \right\}_{x=L} =g&\\
  \end{array}
\end{eqnarray}

where $m=\rho_1h_1+ \rho_3 h_3,$ $K_1=\frac{\rho_1 h_1^3}{12} +\frac{\rho_3 h_3^3}{12},$ and $K_2=\frac{\alpha^1 h_1^3}{12}+\frac{\alpha^3 h_3^3}{12}.$
\section{Well-posedness of the fully dynamic model}

\label{Sec-IV}

In this section, we consider the existence and uniqueness of solutions to (\ref{reduced-nog},\ref{reduced-BC-nog}). First we eliminate the variable $\varphi^1$ by solving the elliptic equation (\ref{reduced-g}) with the associated boundary conditions. Defining $D_x^2 \phi = \phi_{xx}$ and its domain
$${\rm Dom} (D_x^2) = \{ \phi \in H^2 (0,L) ,\quad \phi_x (0)=\phi_x (L)=0 \}, $$
and the operator $P_\xi$ by
\begin{eqnarray}\label{Lgamma}P_{\xi}:=\left(-\xi D_x^2+I\right)^{-1}.\end{eqnarray}
It is well-known that ${P_\xi}$ is a non-negative and a compact operator on $\cLtwo$ \cite{O-M3}. Therefore, the elliptic equation in (\ref{reduced-g}) has the solution
\begin{eqnarray}\label{solutio}
\phi^1=\left\{
  \begin{array}{ll}
   \frac{{\gamma}}{{\ep3}}~P_{\xi} v^3_x , & \sigma_s(t)\equiv 0, i_s^1(x,t)\ne 0, \\
   \frac{{\gamma}}{{\ep3}}~P_{\xi} v^3_x + \frac{\sigma_s(t)}{\ep3 h_3} \chi_{[0,L]}(x)+K , & \sigma_s(t)\ne 0, i_s^1(x,t)\equiv 0.~~
  \end{array}
\right.~~~
 \end{eqnarray}
where $\chi_{[0,L]}(x)$ is the characteristic equation of the interval $[0,L],$ and  $K$ is an arbitrary constant ($\dot K=0$ due to \ref{elec-cont2} and Lemma \ref{pxi}).
Using (\ref{solutio}), the system (\ref{reduced-nog})-(\ref{reduced-BC-nog}) is simplified to
\begin{eqnarray}
\label{reduced-g} \left\{
  \begin{array}{ll}
    \rho_1h_1 \ddot v^1  -\alpha^1 h_1 v^1_{xx}  - G_2  \phi^2 =  f^1,   & \\
  \rho_3h_3 \ddot v^3       -\alpha^3 h_3 v^3_{xx}  + G_2 \phi^2 -{\gamma}h_3 \left(\frac{{\gamma}}{{\ep3}}~(P_{\xi} v^3_x)_x +  \dot \eta_x \right) &\\
  \quad\quad\quad\quad=  \frac{\gamma\sigma_s(t)}{\ep3 } \left(\delta(x)-\delta(x-L)\right) + f^3,   & \\
  m \ddot w - K_1 \ddot{w}_{xx}+K_2w_{xxxx} -  H G_2   \phi^2_x=f&\\
    \xi h_3 \ep3 \ddot \theta + \mu h_3 \left(\theta -\eta_x\right)+ \gamma \xi h_3 (P_{\xi} \dot v^3_x)_x=i_s^1(t) & \\
   \ep3 h_3 \ddot \eta + \mu h_3  \left(\theta-\eta_x\right)_x + \gamma h_3  \left( P_{\xi} -I\right) \dot v^3_x=0&\\
     \phi^2 =\frac{1}{h_2}\left(-v^1+v^3\right) + \frac{H}{h_2}w_x.&
  \end{array} \right.
  \end{eqnarray}
with the initial and  boundary  conditions
\begin{eqnarray}
&&  \begin{array}{ll}
\nonumber   \left\{v^1, v^3, w, w_x \right\}_{x=0}=0, ~~   \left\{\alpha^1 h_1v^1_x\right\}_{x=L}=g^1, &\\
  \end{array} \\
&&  \begin{array}{ll}
\nonumber \left\{ \alpha^3 h_3 v^3_x  +{\gamma}h_3\left( \frac{{\gamma}}{{\ep3}}~P_{\xi} v^3_x+\dot\eta\right)\right\}_{x=L}=g^3, ~~\left\{  \theta, ~\eta_x \right\}_{x=0,L}=0&\\
  \end{array} \\
  &&  \begin{array}{ll}
\nonumber  \left\{K_2 w_{xx}  \right\}_{x=L} =-M,~~\left\{K_1 \ddot w -K_2 w_{xxx} + G_2 H \phi^2 \right\}_{x=L} =g&\\
  \end{array} \\
&&  \begin{array}{ll}
\label{reduced-BC-g}  (v^1,v^3, w, \theta, \eta, \dot v^1, \dot v^3, \dot w, \dot \theta, \dot \eta)(x,0) =(v^1_0, v^3_0, w_0, \theta_0, \eta_0, v^1_1, v^3_1, w_1, \theta_1,\eta_1).&
  \end{array}
\end{eqnarray}
Note that the external forces $g^3$ and $\sigma_s$ control the longitudinal vibrations in the third layer in a similar fashion (too many controllers). At this moment we assume that $g_3\equiv 0$ since the third layer is a piezoelectric layer and it is more appropriate to control this layer by an electrical external force.

Define $\mX=\Ltwo(0,L),$
\begin{eqnarray}
\nonumber && H^1_L(0,L)=\{f\in H^1(0,L)~:~ f(0)=0\}, \\
\nonumber &&H^2_L(0,L)=\{f\in H^2(0,L)\cap H^1_L(0,L)~:~f_x(L)=0 \},\\
 \nonumber   && \mathrm V =  \left\{ \y \in   (H^1_L(0,L))^2\times {H}^2_L(0,L) \times H^1_0(0,L) \times H^1(0,L)\right. \\
\label{H-cur} && \quad\quad\left. -\xi(y_4)_x + y_5= -\xi (y_9)_x + y_{10}=0\right\}
    \end{eqnarray}
and  the complex linear spaces
\begin{eqnarray}
\nonumber  \mc{H} = \mathrm V \times \mathrm H, \quad \mathrm H&=&(\cLtwo)^2 \times  H^1_L(0,L) \times   (\cLtwo)^2.
\end{eqnarray}

The natural energy associated with (\ref{reduced-g}) is
\begin{eqnarray}
  \nonumber &&\mathrm{E}(t)=\frac{1}{2}\int_0^L \left\{\rho_1 h_1  |\dot v^1|^2 + \rho_3 h_3|\dot v^3|^2 + m |\dot w|^2 + K_1 |\dot w_x|^2 + \xi h_3 \ep3 |\dot\theta|^2 \right.\\
  \nonumber && + h_3 \ep3  |\dot \eta|^2 +  {\alpha}^1 h_1 |v^1_x|^2  + \alpha^3 h_3 |v^3_x|^2+\frac{{\gamma}^2  }{\ep3}({P_\xi} v^3_x) v^3_x + K_2 |w_{xx}|^2 \\
\label{Energy-class} &&\left. + \frac{G_2}{h_2}  |-v^1+v^3+H w_x|^2+ \mu h_3  | \theta -\eta_x|^2 \right\} dx.
\end{eqnarray}
This motivates the  definition of the inner product on $\mc H:$
{ \small{\begin{eqnarray}
\nonumber && \left<\left[ \begin{array}{l}
 u_1 \\
 \vdots \\
 u_{10}\\
  \end{array} \right], \left[ \begin{array}{l}
 v_{1} \\
 ~\vdots \\
 v_{10}\\
  \end{array} \right]\right>_{\mc H}= \left<\left[ \begin{array}{l}
 u_6\\
 ~\vdots\\
  u_{10}
 \end{array} \right], \left[ \begin{array}{l}
 v_6\\
 ~\vdots\\
 v_{10}
 \end{array} \right]\right>_{\mathrm H}+ \left<\left[ \begin{array}{l}
 u_1 \\
 ~\vdots \\
  u_5
 \end{array} \right], \left[ \begin{array}{l}
 v_1 \\
 ~\vdots\\
 v_5
 \end{array} \right]\right>_{\mathrm V}\\
\nonumber && =\int_0^L \left\{\rho_1 h_1  u_6 \bar v_6 + \rho_3 h_3 u_7 \bar v_7 + m u_8 \bar v_8  + K_1 (u_8)_x (\bar v_8)_x  + \xi h_3 \ep3 u_9\bar v_9\right.\\
\nonumber  &&   + \ep3 h_3 u_{10} \bar v_{10}  +  {\alpha}^1 h_1 (u_1)_x (\bar v_1)_x + \alpha^3 h_3 (u_2)_x (\bar v_2)_x+\frac{{\gamma}^2  }{\ep3}({P_\xi} (u_2)_x) (\bar v_2)_x \\
\nonumber && + K_2 (u_3)_{xx} (\bar v_3)_{xx}+ \frac{G_2}{h_2}  (-u_1+u_2+H (u_3)_x) (-\bar v_1+\bar v_2+H (\bar v_3)_x)\\
\label{inner-pro}&& \left.+ \mu h_3 ( u_4 -(u_5)_x) (\bar v_4-(\bar v_5)_x)  \right\} dx.~~~~\quad\quad
 \end{eqnarray}}}
  Note that $P_\xi$ is a nonnegative operator, and  $u_4-(u_5)_x=u_4+\xi (u_4)_{xx}$ and $\bar v_4 - (\bar v_5)_x= \bar v_4 + \xi (\bar v_4)_{xx}$ by the definition of $\mathrm H.$ As well, the term $\frac{G_2}{h_2}  (-u_1+u_2+H (u_3)_x) (-\bar v_1+\bar v_2+H (\bar v_3)_x)$ is coercive and continuous in $\mathrm H$, i.e. see \cite{Hansen3}. Therefore (\ref{inner-pro}) is a valid  inner product on $\mathrm H$ and so defines a norm.
 It is straightforward to verify that $E$  as defined in (\ref{Energy-class}) is the norm induced by (\ref{inner-pro}). It can also easily be shown that $\mathrm H$ with this norm is complete. This follows from the fact that the gauge constraints in $\mathrm H$ are satisfied weakly \cite{O-M3}.

Let $\vec y=(v^1, v^3, w, \theta, \eta)$ be the smooth solution of the system of (\ref{reduced-g})-(\ref{reduced-BC-g}).  Multiplying the equations in (\ref{reduced-g}) by $\tilde y_1, \tilde y_2 \in H^1_L(0,L), ~ \tilde y_3 \in H^2_L(0,L),$ $\tilde y_4\in H^1_0(0,L),$ and $ \tilde y_5\in H^1(0,L),$ respectively, where $\xi (\tilde y_4)_x=\tilde y_5,$ and integrating by parts yields
\begin{eqnarray}
\label{reduced-gg}  \begin{array}{ll}
   \int_0^L \left\{ \rho_1h_1 \ddot v^1 \tilde y_1 +\alpha^1 h_1 v^1_{x} (\tilde y_1)_x  - G_2  \phi^2\tilde y_1\right\}~dx -g_1(t)\tilde y_1(L)= 0,   & \\
  \int_0^L \left\{\rho_3h_3 \ddot v^3  \tilde y_2     +\alpha^3 h_3 v^3_{x} (\tilde y_2)_x  + G_2 \phi^2 \tilde y_2 +{\gamma}h_3 \left(\frac{{\gamma}}{{\ep3}}~(P_{\xi} v^3_x) +  \dot \eta \right) (\tilde y_2)_x\right\}~dx &\\
  \quad\quad\quad\quad=  -\frac{\gamma\sigma_s(t)}{\ep3 } \tilde y_2(L),   & \\
 \int_0^L \left\{ m \ddot w \tilde y_3+ K_1 \ddot{w}_{x} (\tilde y_3)_{x}+K_2w_{xx} (\tilde y_3)_{xx} +  H G_2   \phi^2 (\tilde y_3)_x\right\}~dx &\\
 \quad\quad\quad\quad-g(t) \tilde y_3(L) + K_2 M(t)(\tilde y_3)_x(L) =0& \\
  \int_0^L \left\{ \xi \ep3 h_3\ddot \theta \tilde y_4+\mu h_3 \left(\theta -\eta_x\right) \tilde y_4- \gamma \xi h_3 (P_{\xi} \dot v^3_x) (\tilde y_4)_x\right\} ~dx=\int_0^L i_s^1(t) \tilde y_4 ~dx & \\
   \int_0^L \left\{ \ep3 h_3 \ddot \eta \tilde y_5 - \mu h_3 \left(\theta-\eta_x\right) (\tilde y_5)_x + \gamma h_3  \left( P_{\xi} -I\right) \dot v^3_x ~\tilde y_5\right\}~dx =0&
  \end{array}
  \end{eqnarray}
Note that, the term $\delta(x)\sigma_s(t)$ in (\ref{reduced-g}) is lost in the weak formulation due to the boundary condition $\tilde y_2(0)=0.$ This explains that the charge actuation only forces boundary condition at $x=L$ as in the case of voltage actuation \cite{Ozkan2}.

Let $\mc M : H^1_L (0,L) \to (H^1_L(0,L))'$ and $D:H^1_L (0,L) \to (H^1_L(0,L))'$ be a linear operators defined by
\begin{eqnarray}\nonumber &&\left< \mc M \psi, \tilde \psi \right>_{(H^1_L(0,L))', H^1_L(0,L)}= \int_0^L \left(m \psi \bar{\tilde \psi} + K_1 \psi_x {\bar{\tilde \psi}}_x\right) dx.\\
\label{def-mcM}&&\left< \mc D \psi, \tilde \psi \right>_{V',V}= \gamma h_3\int_0^L \left( (\psi_5)_x(\bar{\tilde \psi}_2)_x- (\psi_2)_x (\bar{\tilde \psi}_5)_x  \right) dx.
\end{eqnarray}
For all $\psi=(\psi_1, \psi_2,\psi_3, \psi_4, \psi_5)^{\rm T}, \tilde \psi= (\tilde \psi_1, \tilde \psi_2, \tilde \psi_3, \tilde \psi_4, \tilde \psi_5)^{\rm T} \in V,$  define the linear operator $A : V \to V'$ by
$\left<A \psi, \tilde\psi\right>_{V',V}= \left<\psi, \tilde \psi\right>_{V}.$
From the Lax-Milgram theorem $\mc M$ and $A$ are  canonical  isomorphisms from $H^1_L(0,L)$ onto  $(H^1_L(0,L))'$ and from $V$ onto $V',$ respectively.

Next we introduce the linear unbounded operator by
$\mc A: {\text{Dom}}(A)\times V\subset \mc H \to \mc H$
where $\mc A= \left[ {\begin{array}{*{20}c}
   O_{5\times 5}  & I_{5\times 5} \\
       M^{-1}A  &  M^{-1}D   \\
\end{array}} \right],$ and
\begin{eqnarray}  \label{defA}{\rm {Dom}}(\mc A) &=&  \{(\vec z, \vec {\tilde z}) \in V\times V, A\vec z + D\dot z\in \mathrm V' \}.
\end{eqnarray}

 Let ${\rm Dom}(\mc A)'$ is the dual of ${\rm Dom}(\mc A)$ pivoted with respect to $\mc H.$ Define the control operator $B,$ and $ B_0 \in \mathcal{L}(\mathbb{C}, \mathrm V'),$
 \begin{eqnarray}
\nonumber  ~ \text{with} ~ B_0 =   \left[ \begin{array}{cccccccc}
\delta(x-L) & I & 0 & 0 & 0 & 0 & 0 & 0 \\
0 & 0 & \frac{-\gamma}{\ep3}\delta(x-L) & I & 0 & 0 & 0 & 0\\
0 & 0 & 0 & 0 & \delta(x-L)& K_2 \delta '(x-L) &  I & 0 \\
0 & 0 & 0 & 0 & 0 & 0 & 0 & \frac{1}{h_3\ep3} I\\
0 & 0 & 0 & 0 & 0 & 0 & 0 & 0
 \end{array} \right],&
 \end{eqnarray}
 \begin{eqnarray}
 \label{defb_0} &\quad B \in \mathcal{L}(\mathbb{C} , ({\rm Dom}(\mc A))'), ~ \text{with} ~ B=   \left[ \begin{array}{c} 0_{5\times 8} \\ B_0 \end{array} \right]&
 \end{eqnarray}
  We have  $\mathrm V'$ which is the dual of $\mathrm V$ pivoted with respect $\mathrm H.$ The dual operators of $B_0$ and $B$ are $B_0^* \in \mc L( \mathrm V,  \mathbb{C})$ and
\begin{eqnarray*}
 B_0^* \psi &=& \left[\psi_1(L), \int_0^L \psi_1(x)dx, \frac{-\gamma}{\ep3}\psi_2(L), \int_0^L \psi_2(x)dx,  \psi_3(L),  -K_2 \psi_3'(L), \right.\\
   &&\left.\int_0^L \psi_3 (x)dx, \frac{1}{h_3\ep3}\int_0^L \psi_4 (x)dx \right]^{\rm T}
\end{eqnarray*}
with $B^*\Phi=(0_{8\times 5}\quad B_0^*)^{\rm T}\Phi.$ Let $M={\rm diag}\left(\rho_1 h_1 I ~~\rho_3 h_3 I~~ \mc M~~\xi  \ep3 h_3 I~~ \ep3 h_3 I\right).$ $M$ is clearly is an  isomorphism from $\mathrm H$ onto $\mathrm H'.$

Writing $\varphi=[v^1, v^3, w, \theta, \eta, \dot v^1,  \dot v^3, \dot w, \dot \theta, \dot\eta]^{\rm T},$ the control  system (\ref{reduced-g})-(\ref{reduced-BC-g})  with all controllers $F=[g^1, g^1, \sigma_s, f^3, g, M, f, i_s^1]^{\rm T} $ can be put into the  state-space form
\begin{subequations}
\label{Semigroup}
\begin{empheq}[left={\phantomword[l]{}{ }  \empheqlbrace}]{align}
&\dot \varphi +  \underbrace{ \left[ {\begin{array}{*{20}c}
   0 & -I_{5\times 5}  \\
    M^{-1} A & M^{-1}D  \\
\end{array}} \right]}_{\mc A}  \varphi =\underbrace{ \left[ \begin{array}{c} 0_{5\times 8} \\ M^{-1}B_0 \end{array} \right]}_B F(t) , &\\
&\varphi(x,0) =  \varphi ^0.
\end{empheq}
\end{subequations}
Note that the term $-\gamma \xi h_3 (P_{\xi} \dot v^3_x) (\tilde y_4)_x$ in (\ref{reduced-gg}) is the same as $=- \gamma h_3 (P_{\xi} \dot v^3_x) \tilde y_5$ in the last equation due to the Coulomb gauge (\ref{Colombo-trans}).

\begin{lem} \label{skew}The operator $\mc{A}$  defined by (\ref{defA}) is maximal monotone in the energy space $\mc H,$ i.e.,
$$\left< \mc A \vec z, \vec z \right> \ge 0, ~~~{\rm for~ all} ~ \vec z\in {\rm Dom}(\mc A),$$
 and ${\rm Range}(I+\mc A)=\mc H.$
\end{lem}

\textbf{Proof:} Let $\vec z \in {\rm Dom}(\mc A).$ A simple calculation using integration by parts and the boundary conditions yields
{\small{\begin{eqnarray}
\nonumber \left<\mc A  \left[ \begin{array}{c}
\vec z_1 \\
\vec z_2 \end{array} \right], \left[ \begin{array}{c}
\vec z_1 \\
\vec z_2 \end{array} \right]\right>_{\mc H\times \mc H} &=& \left<  \left[ \begin{array}{c}
-\vec z_2 \\
M^{-1} A \vec z_1 + M^{-1}Dz_2 \end{array} \right], \left[ \begin{array}{c}
\vec z_1 \\
\vec z_2 \end{array} \right]\right>_{\mc H}\\
\nonumber  &=& \left<-\vec z_2, \vec z_1\right>_{V\times V}  + \left<M^{-1}A\vec z_1 +  M^{-1} Dz_2, \vec z_2\right>_{H\times H}\\
\nonumber &=& -\overline{\left<A \vec z_1, \vec z_2\right>_{V'\times V}}  + \left<A\vec z_1 + Dz_2, \vec z_2\right>_{H'\times H}.
\end{eqnarray}
}}
Since $\vec z=\left[ \begin{array}{c}
\vec z_1 \\
\vec z_2 \end{array} \right] \in {\rm Dom}(\mc A) $, then $A\vec z _1 + D z_2 \in V'$ and $\vec z_2 \in V$ so that
$$\left<A\vec z_1 + Dz_2, \vec z_2\right>_{H'\times H}=\left<A\vec z_1+Dz_2 , \vec z_2\right>_{V'\times V}.$$
Since ${\rm Re}\left<Dz_2 , \vec z_2\right>_{V'\times V}=0$ by (\ref{def-mcM}),
${\rm Re}\left<\mc A  \vec z, \vec z\right>_{\mc H\times \mc H}=0.$

We next verify the range condition. Let $\vec z=\left[ \begin{array}{c}
\vec z_1 \\
\vec z_2 \end{array} \right]\in \mc H.$ We prove that there exists a $\vec y=\left[ \begin{array}{c}
\vec y_1 \\
\vec y_2 \end{array} \right]\in \mc {\rm Dom} (\mc A )$ such that $(I+\mc A) \vec y=\vec z.$
A simple computation shows that proving this is equivalent to proving ${\rm Range} (M+A+D)=H',$ i.e., for every $\vec f \in H'$   there exists a unique solution $\vec z \in H$ such that $(M+A+D)\vec z=\vec f.$
This obviously follows from the Lax Milgram's theorem.

\begin{prop}\label{prop-dom}Let  $\vec u=(\vec y, \vec z)\in \mc H.$ Then $\vec u \in {\rm Dom}(\mc A)$  if and only if the following conditions hold:
\begin{eqnarray}
\nonumber && \vec y\in (H^2(0,L)\cap H^1_L(0,L))^2 \times (H^3(0,L) \cap H^2_L(0,L))\\
\nonumber &&\quad\quad\times (H^2(0,L)\cap H^1_0(0,L)) \times H^2(0,L),~~{\rm and}~~\vec z\in V,\\
\nonumber && \xi (y_4)_x=y_5, ~ \xi (z_4)_x = z_5~~{\rm such ~~ that}~~ (y_5)_x(0)=(y_5)_x(L)=0,~~{\rm and }\\
\label{Dom} && (y_1)_{x}=\alpha^3 h_3 (y_2)_x  +{\gamma}h_3\left( \frac{{\gamma}}{{\ep3}}~P_{\xi} (y_2)_x+z_5\right)=(y_3)_{xx}\left. \right|_{x=L}=0.
\end{eqnarray}
Moreover, the resolvent of $\mc A$ is compact in the energy space $\mc H.$
\end{prop}

  \textbf{Proof:} Let $\vec {\tilde u}= \left( \begin{array}{c}
\vec {\tilde y} \\
\vec {\tilde z}\\
 \end{array} \right)\in \mc H$ and $\vec {u}= \left( \begin{array}{c}
\vec {y} \\
\vec {z}\\
 \end{array} \right)\in {\rm Dom}(\mc A)$ such that $\mc A \vec u=\vec {\tilde u}.$ Then we have
 \begin{eqnarray}
 \nonumber -\vec z&=& \vec {\tilde y} \in V,\\
 \label{dumb}   A\vec y +D\vec z &=&M\vec{\tilde z},
 \end{eqnarray}
  and therefore,
  \begin{eqnarray}
\label{dumb1} \left< \vec y, \vec \varphi\right>_{V}+\left<Dz,\varphi\right>_{\mathrm V',V}= \left< \vec {\tilde z}, \vec \varphi \right>_{H} ~~{\rm for ~~all ~~} \vec \varphi \in V.
\end{eqnarray}
Let  $\vec \psi=[\psi_1,\ldots,\psi_5]^{\rm T} \in (C_0^\infty(0,L))^4.$  We define $\varphi_i=\psi_i$ for $i=1,2,4,5,$ and $\varphi_3=\int_0^x \psi_3 (s)ds.$ Since $\vec \varphi\in V,$ inserting $\vec \varphi$ into (\ref{reduced-gg}) and using the  gauge condition  yield
  \begin{footnotesize}
\begin{eqnarray}
\label{reduced-what}
  \begin{array}{ll}
   \int_0^L \left\{ -\alpha^1 h_1 v^1_{xx}  \bar\psi_1-\alpha^3 h_3 v^3_{xx} \bar \psi_2 - \frac{{(\gamma^3)^2 h_3}}{{\ep3}}~(P_{\xi} v^3_x)_x   \bar \psi_2- K_2w_{xxx} \bar \psi_3 + G_2  h_2\phi^2 \bar \psi^2 \right.,   & \\
    \left.  +\mu h_3 \left(\theta -\xi \theta_{xx}\right) (\bar \psi_4) + \mu h_3 \left(\frac{1}{\xi}\eta -\eta_{xx}\right) \bar \psi_5 \right\}~dx =\int_0^L\left\{  \left(\int_1^x m \tilde z_3  ds + K_1 (\tilde z_3)_x\right) \bar \psi_3  \right.&\\
    \left. \rho_1 h_1 z_1 \bar \psi_1+ \rho_3h_3 z_2 \bar \psi_2 + \xi \ep3 h_3 z_4 \bar \psi_4 + \ep3 h_3 z_5 \bar \psi_5-\gamma h_3 (z_7)_x\bar \psi_5 - \gamma h_3 (z_5)_x \bar \psi_2   \right\}~dx=0 &
  \end{array}
  \end{eqnarray}
\end{footnotesize}
for all $\vec \psi\in (C_0^\infty(0,L))^4.$
Therefore it follows that $\vec y\in (H^2(0,L)\cap H^1_L(0,L))^2 \times (H^3(0,L) \cap H^2_L(0,L))\times (H^2(0,L)\cap H^1_L(0,L))^2.$
Next let $\vec \psi \in \mathrm H.$  We define
\begin{eqnarray} \label{dumber}\varphi_i=\int_0^x \psi_i(s)ds,\quad i=1,\ldots,4.
\end{eqnarray}
Obviously $\vec \varphi \in \rm V. $ Then plugging (\ref{dumber}) into (\ref{dumb1}) yields
\begin{eqnarray} \nonumber  &&\alpha^1 h_1  (y_1)_{x}(L) \bar \psi_1(L) + \left(\alpha^3 h_3 (y_2)_x + \frac{(\gamma^3)^2 h_3}{\ep3} (P_\xi)(y_2)_x (L) + \gamma^3 h_3z_5(L)\right) \bar \psi_2 (L) \\
\nonumber &&\quad+ K_2 (y_3)_{xx}(L) \bar \psi_3(L) +\y_4 (L) \bar \psi_4(L) +  (y_5)_x(L)  \bar \psi_5(L)=0
\end{eqnarray}
for all $\psi\in \mathrm H.$ Hence,
$$(y_4)_{xx} (L)=(y_1)_x(L)=(y_2)_x(L)=(y_3)_x(L)=0.$$

We have the following well-posedness theorem for (\ref{Semigroup}).
\begin{thm}\label{w-pf}
Let $T>0,$ and $F\in (\Ltwo(0,T))^{8}.$ For any $\varphi^0 \in \mc H,$ $\varphi\in C[[0,T]; \mc H],$ and there exists a positive constants $c_1(T)$
such that (\ref{Semigroup}) satisfies
      \begin{eqnarray}\label{conc}\|\varphi (T) \|^2_{\mc H} &\le& c_1 (T)\left\{\|\varphi^0\|^2_{\mc H} + \|F\|^2_{(\Ltwo(0,T))^8}\right\}.
      \end{eqnarray}
\end{thm}
\vspace{0.1in}

\textbf{Proof:} The operator $\mc A: {\rm {Dom}}(\mc A) \to \mc H$ is a maximal monotone operator by Lemma \ref{skew} and Proposition \ref{prop-dom}. Therefore it is an  infinitesimal generator of $C_0-$semigroup of contractions by Hille-Yosida theorem. Since $$ \left<M^{-1}B F, \tilde \psi\right>_{V,V}=\left<B F, \tilde \psi\right>_{V',V}<\infty$$ by the Sobolev's trace inequality, $B$ is an admissible control operator for the semigroup $\{e^{\mc A t}\}_{t\ge 0}$ corresponding to (\ref{Semigroup}), and hence the conclusion (\ref{conc})  follows. $\square$

\begin{rmk} Note that the the underlying control operators in each actuation differ substantially. In the case of current actuation, the control operator is bounded in $\mathrm H$ whereas in the case of charge actuation it is unbounded. Therefore the charge actuation case is more similar to the voltage actuation case in this manner.
\end{rmk}
\noindent{\bf Quasi-static model:} The quasi-static assumes that $\dot A=0$ but $A$ may be nonzero, i.e. $E=-\nabla \phi,$ and $\theta,\eta$ may be nonzero. The Coulomb gauge  (\ref{Colombo-trans})is still valid. Running the Hamilton's principle by taking this into account together with the Rao-Nakra assumptions, i.e. $\phi_2, \alpha^2\to 0$ and  using the definition of $\xi,$ yields the following initial-boundary value problem:
\begin{eqnarray}
\label{quasi} \left\{
  \begin{array}{ll}
    \rho_1h_1 \ddot v^1  -\alpha^1 h_1 v^1_{xx}  - G_2  \phi^2 = 0,   & \\
  \rho_3h_3 \ddot v^3       -\alpha^3 h_3 v^3_{xx}  + G_2 \phi^2 - \frac{{\gamma^2}}{{\ep3 h_3}}~(P_{\xi} v^3_x)_x =  -\frac{\gamma\sigma_s(t)}{\ep3 } \delta(x-L) ,   & \\
  m \ddot w - K_1 \ddot{w}_{xx}+K_2w_{xxxx} -  H G_2   \phi^2_x=0&\\
      \mu h_3 \left(\theta -\eta_x\right)+ \gamma \xi h_3 (P_{\xi} \dot v^3_x)_x=i_s^1(t) & \\
   \mu h_3  \left(\theta-\eta_x\right)_x + \gamma h_3  \left( P_{\xi} -I\right) \dot v^3_x=0&\\
     \phi^2 =\frac{1}{h_2}\left(-v^1+v^3\right) + \frac{H}{h_2}w_x.&
  \end{array} \right.
  \end{eqnarray}
with the initial and boundary  conditions
\begin{eqnarray}
&&  \begin{array}{ll}
\nonumber   \left\{v^1, v^3, w, w_x \right\}_{x=0}=0, ~~   \left\{\alpha^1 h_1v^1_x\right\}_{x=L}=g^1,&\\
  \end{array} \\
  &&  \begin{array}{ll}
\nonumber \left\{ \alpha^3 h_3 v^3_x  + \frac{{\gamma^2}}{{\ep3 h_3}}~P_{\xi} v^3_x\right\}_{x=L}=0, ~~\left\{  \theta,~  \eta_x \right\}_{x=0,L}=0, &\\
  \end{array} \\
      &&  \begin{array}{ll}
\nonumber      \left\{K_2 w_{xx}  \right\}_{x=L} =-M,\quad \left\{K_1 \ddot w -K_2 w_{xxx} + G_2 H \phi^2 \right\}_{x=L} =g&\\
  \end{array} \\
&&  \begin{array}{ll}
\label{quasi-BC}  (v^1,v^3, w, \theta, \eta, \dot v^1, \dot v^3, \dot w)(x,0) =(v^1_0, v^3_0, w_0, \theta_0, \eta_0, v^1_1, v^3_1, w_1).&
  \end{array}
\end{eqnarray}

 \begin{lem} \cite[Lemma 2.1]{O-M3} \label{pxi} Let ${\rm Dom} (D_x^2)=\{w\in H^2(0,L) ~:~ w_x(0)=w_x(L)=0\}.$ Define the operator $J=\frac{1}{\xi}(P_\xi-I). $ Then $J$ is continuous, self-adjoint and and non-positive on $\cLtwo.$ Moreover, for all $w\in {\rm Dom} (P_\xi), $
 $$J w = D_x^2 ~P_\xi= D_x^2 (I-\xi D_x^2)^{-1}w.$$
\end{lem}

Note that the last equation can easily be derived as the first derivative of the $\theta-$equation is taken by using Lemma \ref{pxi}, the surface continuity condition (\ref{elec-cont2}), and the gauge condition (\ref{Colombo-trans}) in a weak sense. This way the last two equations can be also written as second order elliptic equations as
\begin{footnotesize}
\begin{eqnarray}
\label{quasi-eliptic} \left\{
  \begin{array}{ll}
      \mu h_3 \left(\theta -\xi \theta_{xx}\right)+ \gamma \xi h_3 (P_{\xi} \dot v^3_x)_x=i_s^1(t) & \\
   \mu h_3   \left(\eta-\xi\eta_{xx}\right) + \gamma \xi h_3 \left( P_{\xi} -I\right) \dot v^3_x=0&\\
  \end{array} \right.
  \end{eqnarray}
\end{footnotesize}

Since the bending and stretching equations is not coupled to $\theta$ or $\eta$ equation directly, this case should be handled delicately in the case of current actuation. The well-posedness of this model can be shown analogously.
\vspace{0.1in}

\noindent {\bf Electrostatic model:} The electrostatic assumption fully ignores magnetic effects, i.e. $A=\dot A=0.$ Note that the surface current density $i_s(t)$ can not survive in the fully dynamic approach  in (\ref{work}). For this reason, current surface density may have to come into the play through the circuit equations. However, the surface charge density $\sigma_s(t)$ survives in (\ref{work}). If we run the Hamilton's principle by taking this into account together with the reduced Rao-Nakra assumptions, i.e. $\phi_2, \alpha^2\to 0,$ and  using the definition of $\xi,$ we obtain the following initial-boundary value problem:
\begin{eqnarray}
\label{electro} \left\{
  \begin{array}{ll}
    \rho_1h_1 \ddot v^1  -\alpha^1 h_1 v^1_{xx}  - G_2  \phi^2 =0,   & \\
  \rho_3h_3 \ddot v^3       -\alpha^3 h_3 v^3_{xx} -\frac{{\gamma^2}}{{\ep3 h_3}}~(P_{\xi} v^3_x)_x  + G_2 \phi^2  = -\frac{\gamma\sigma_s(t)}{\ep3 } \delta(x-L),   & \\
  m \ddot w - K_1 \ddot{w}_{xx}+K_2w_{xxxx} -  H G_2   \phi^2_x=0&\\
     \phi^2 =\frac{1}{h_2}\left(-v^1+v^3\right) + \frac{H}{h_2}w_x.&
  \end{array} \right.
  \end{eqnarray}
with the natural boundary  conditions at $x=0,L$ and initial conditions
\begin{eqnarray}
&&  \begin{array}{ll}
\nonumber   \left\{v^1, v^3, w, w_x  \right\}_{x=0}=0, ~  \left\{\alpha^1 h_1v^1_x\right\}_{x=L}=g^1,&\\
  \end{array} \\
  &&  \begin{array}{ll}
\nonumber  \left\{ \alpha^3 h_3 v^3_x  +\frac{{\gamma^2}}{{\ep3 h_3}}~P_{\xi} v^3_x\right\}_{x=L}=0,&\\
  \end{array} \\
         &&  \begin{array}{ll}
\nonumber  \left\{K_2 w_{xx}  \right\}_{x=L} =-M, ~  \left\{K_1 \ddot w -K_2 w_{xxx} + G_2 H \phi^2 \right\}_{x=L} =g&\\
  \end{array} \\
&&  \begin{array}{ll}
\label{electro-BC}  (v^1,v^3, w, \dot v^1, \dot v^3, \dot w)(x,0) =(v^1_0, v^3_0, w_0, v^1_1, v^3_1, w_1).&
  \end{array}
\end{eqnarray}
Note that this case is not much different from the voltage-actuation case \cite{Ozkan2,Ozkan3}. The difference arises from the term involving the $P_\xi$ term in (\ref{electro}). This is only because the voltage actuated model ignores the induced effect of the electric field, i.e. see (\ref{scalarpot}). We account for this effect by choosing quadratic-through thickness assumption. These effects, in fact, turns out to make the piezoelectric beam more stiff since $P_\xi$ is a positive operator. This was first observed in \cite{Hansen}. The well-posedness of this model can be shown analogously.

  \section{Future Research}

  \begin{figure}
\begin{center}
\includegraphics[height=4cm]{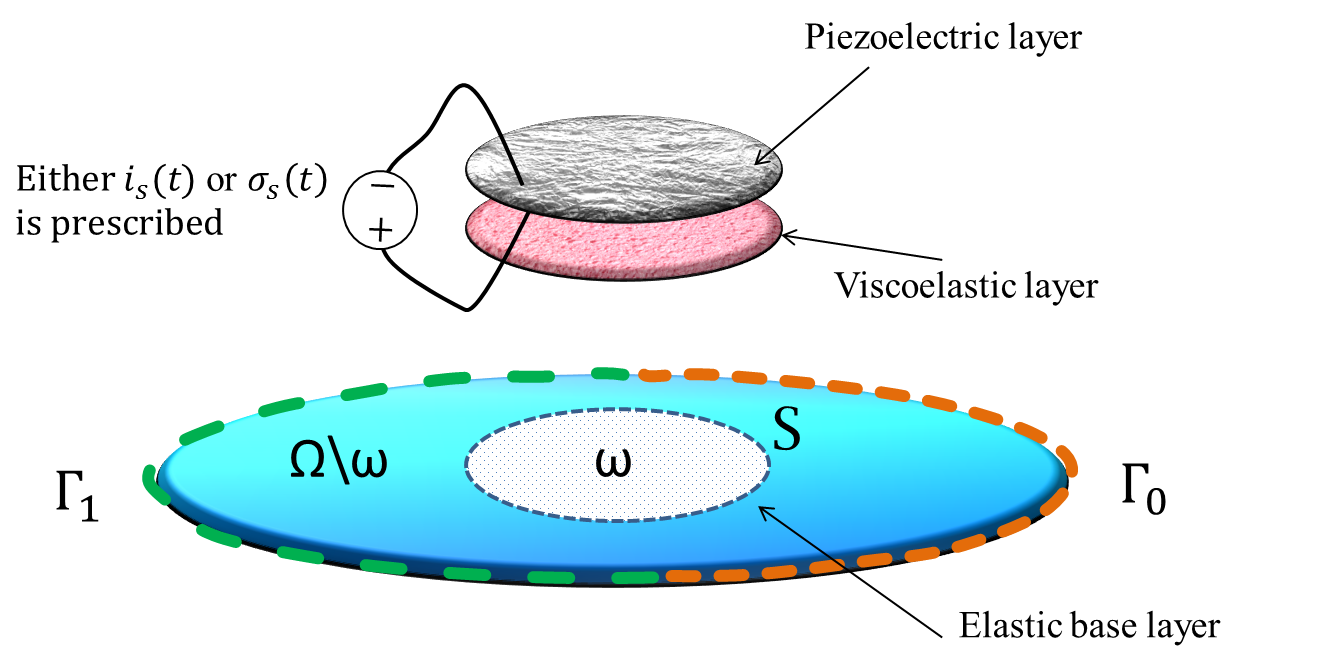}    
\caption{ACL-patch on an elastic host layer.}  
\label{ACL-patch}                                 
\end{center}                                 
\end{figure}
A relevant problem is  the ACL-patch problem on a stiff elastic layer, see Figure \label{ACL-patch}. The model corresponding to this problem is obtained in a similar fashion by following the methodologies in this paper and \cite{accpaper}. However, the exact controllability and stabilizability problem is a bit more challenging relying on the location of the patch and the material parameters. This is currently under investigation.

The three-layer ACL beam model obtained in this paper can be analogously extended to the multilayer ACL beam model with the same assumptions as in \cite{Hansen3}. The piezoelectric layers can be even actuated by charge, voltage, or current sources simultaneously. This is the subject of the future research.

\bibliographystyle{plain}        

\end{document}